\documentclass[12pt]{article}
\usepackage{amssymb,amsmath,amsthm}
\usepackage[dvips]{graphicx,xcolor}
\usepackage{here}

\newcommand{\pf}{\vspace{3mm}\noindent{\it Proof)}\hspace{5mm}}

\newcommand{\R}{{\bf R}}

\newcommand{\p}{\partial}

\newcommand{\lsm}{\lesssim}
\newcommand{\ep}{\varepsilon}

\newcommand{\al}{\alpha}
\newcommand{\bt}{\beta}

\newcommand{\lm}{\lambda}

\newcommand{\dsp}{\displaystyle}
\newcommand{\lr}{\longrightarrow }

\newcommand{\nb}{\nabla}
\newcommand{\fa}{\frac}
\newcommand{\sr}{\sqrt}
\newcommand{\lp}[1]{\left(#1 \right)}
\newcommand{\lb}[1]{\left\{#1 \right\}}
\newcommand{\lbt}[1]{\left[#1 \right]}

\begin{document}
 
\addtocounter{footnote}{1}
 \title{}
 
\begin{center}
{\large\bf 
Blow up of solutions of semilinear wave equations \\
in  accelerated expanding Friedmann-Lema\^itre-Robertson-Walker spacetime\\

}
\end{center}
\vspace{3mm}
\begin{center}

Kimitoshi Tsutaya$^\dagger$ 
        and Yuta Wakasugi$^\ddagger$ \\
\vspace{1cm}

 $^\dagger$Graduate School  of Science and Technology \\
Hirosaki University  \\
Hirosaki 036-8561, Japan\\
\footnotetext{* The research was supported by JSPS KAKENHI Grant Number JP18K03351. }

\vspace{5mm}
        $^\ddagger$
Graduate School of Engineering \\
Hiroshima University \\
Higashi-Hiroshima, 739-8527, Japan

\end{center}

\vspace{3mm}
\begin{abstract}
Consider a nonlinear wave equation  for a massless scalar field  with self-interaction 
in the spatially flat Friedmann-Lema\^itre-Robertson-Walker spacetimes. 
For the case of accelerated expansion, 
we show that blow-up in a finite time occurs for the equation with arbitrary power nonlinearity as well as upper bounds of the lifespan of blow-up solutions. 
Comparing to the case of the Minkowski spacetime, 
we discuss how the scale factor affects the lifespan of blow-up solutions of the equation. 
\end{abstract}  


\section{Introduction.}
\addtolength{\baselineskip}{2.4mm}
\quad 
This is the third in a series of papers concerned with  
the wave equation describing a mass-less scalar 
field with self-interaction in the Friedmann-Lema\^itre-Robertson-Walker (FLRW) spacetimes. 
The spatially flat FLRW metric is given by 
\[
g: \; ds^2=-dt^2+a(t)^2d\sigma^2, 
\]
where the speed of light is equal to $1$, $d\sigma^2$ is the line element of $n$-dimensional Euclidean space and $a(t)$ is the scale factor, which describes expansion or contraction of the spatial metric. 

If we solve the Einstein equation with the energy-momentum tensor for the perfect fluid under an equation of state, we obtain 
\begin{equation}
a(t)=ct^{\fa 2{n(1+w)}}  
\label{scw} 
\end{equation}
with some constant $c$, where $w$ is a proportional constant in the equation of state. See \cite{CGLY} for details.  The present paper treats the case $-1<w\le 2/n-1$ for $n\ge 2$, which describes an accelerated expanding universe. We note that the case $w=2/n-1$ corresponds to a uniformly accelerated expansion.

For the flat FLRW metric with \eqref{scw},  
the semilinear wave equation $\Box_g u\\
=|g|^{-1/2}\p_\al(|g|^{1/2}g^{\al\bt}\p_\bt)u= -|u|^p$ with $p>1$ becomes 
\begin{equation}
u_{tt}-\fa 1{t^{4/(n(1+w))}}\Delta u+\fa 2{(1+w)t}u_t=|u|^p,   \quad  x\in \R^n 
\label{ore}
\end{equation}
where  $\Delta=\p_1^2+\cdots \p_n^2, \; \p_j=\p/\p x^j, \; j=1,\cdots,n, \; (x^1,\cdots, x^n)\in \R^n$
.  
Our aim of this paper is to show that blow-up in a finite time occurs for the equation above as well as upper bounds of the lifespan of the blow-up solutions. 

We first consider the following Cauchy problem
in order to compare with the related known results including the case of the Minkowski spacetime: 
\begin{equation}
\begin{cases}
\dsp
u_{tt}-\fa 1{t^{2\al}}\Delta u+\fa\mu{t}u_t=|u|^p, &\qquad  t>1, \; x\in \R^n  \\
  & \\
u(1,x)=\ep u_0(x), \; u_t(1,x)=\ep u_1(x), &\qquad x\in \R^n,  
\end{cases}
\label{eq}
\end{equation}
where $\al$ and $\mu$ are nonnegative constants, and $\ep>0$ is a small parameter.  
We then return to \eqref{ore}.

Let $T_\ep$ be the lifespan of solutions of \eqref{eq}, that is, $T_\ep$ is the supremum of $T$ such that  \eqref{eq} have a solution for $x\in \R^n$ and $1\le t<T$. 
Let $p_F(n)$ denote the Fujita exponent $p_F(n)=1+2/n$, and $p_S(n)$ the Strauss exponent which is the positive root of the equation 
\begin{equation}
\gamma_S(n,p)\equiv -(n-1)p^2+(n+1)p+2=0. 
\label{snp}
\end{equation}
For the special case $\al=0$, we refer to \cite{IkSo18} and references therein. 

The authors have recently proved blow-up in a finite time and upper estimates of the lifespan for the case $0\le \al <1$. 
Let 
\begin{align}
\gamma(n,p,\al,\mu)=-p^2\lp{n-1+\fa{\mu-\al}{1-\al}}+p\lp{n+1+\fa{\mu+3\al}{1-\al}}+2
\label{gm} 
\end{align}
and $p_c(n,\al,\mu)$ be the positive root of the equation $\gamma(n,p,\al,\mu)=0$. 
It is shown in \cite{TW1} that 
\begin{align}
&T_\ep\le C\ep^{\fa{-2p(p-1)}{(1-\al)\gamma(n,p,\al,\mu)}} \quad \mbox{ if }\quad 1<p<p_c(n,\al,\mu), 
\label{TW1} \\
&T_\ep\le \exp(C\ep^{-p(p-1)}) \quad \mbox{ if }\quad p=p_c(n,\al,\mu)>p_F(n(1-\al))=1+2/(n(1-\al)). 
\label{TW2}
\end{align}
We remark that if $\mu=\al=0$, the $p_c(n,0,0)$ coincides with the Strauss exponent $p_S(n)$. 
If $\al=2/3$ and $\mu=2$, which is the Einstein-de Sitter spacetime if $n=3$,  
 our blow-up range of $p$ is the same as $1<p<p_{cr}(n)$ in Galstian and Yagdjian \cite{GY}. 
We note that these upper bounds stated so far are sharp 
if the power $p$ is dominated by the Strauss exponent or $p_c(n,\al,\mu)$. 
We can say in this case that the lifespan of solutions is {\it wavelike} as mentioned in \cite{LTW}.  

On the other hand, if the critical exponent for \eqref{eq} is affected by the Fujita exponent, the upper bound of the lifespan is expected to become sharper. 
In fact, we have proved in \cite{TW2} the following upper bounds of the lifespan for $0\le \al<1$:  
\begin{align*}
&T_{\varepsilon} \le C \varepsilon^{- \frac{p-1}{2-n(1-\alpha)(p-1)}}
 \quad \mbox{ if }\quad 1<p<p_F(n(1-\alpha)),\\
&T_{\varepsilon} \le \exp\left( C\varepsilon^{-(p-1)} \right)
  \quad \mbox{ if }\quad p=p_F(n(1-\alpha)).
\end{align*}
We can say that this case corresponds to the {\it heatlike} lifespan as in \cite{LTW}. 
These results generalize ones given in \cite{Wa1,Wa2} and Ikeda et al \cite{ISW}.

In this paper we treat the remaining case $\al\ge 1$, and show upper estimates of the lifespan of the blow-up solutions. 
We then apply our results to the original equation \eqref{ore}. 
Our aim is especially to clarify the difference with the case of the Minkowski spacetime and also how the scale factor affects the lifespan of the solution.

The following theorem is our main result, which states that blow-up in a finite time can occur for all $p>1$. 
  
\noindent
{\bf Theorem 1.1.} 
{\it 
 Let $n\ge 2, \; \al\ge 1, \; \mu\ge 0$ and $p>1$.
Assume that $u_0\in C^2(\R^n)$ and $u_1\in C^1(\R^n)$ are nonnegative, nontrivial and $\mbox{supp }u_0, \mbox{supp }u_1\subset \{|x|\le R\}$ with $R>0$. 
Suppose that the problem \eqref{eq} has  a classical solution $u\in C^2([1,T)\times\R^n)$.  
Then $T<\infty$ and there exists a constant $\ep_0>0$ depending on 
$p,\al,\mu,R,u_0,u_1$ such that $T_\ep$ has to satisfy 
\begin{align*}
&T_\ep^2(\ln T_\ep)^{-n(p-1)}\le C\ep^{-(p-1)} &&\mbox{if }\al=1, \\
&T_\ep\le C\ep^{-(p-1)/2} &&\mbox{if }\al>1
\end{align*}
for $0<\ep\le \ep_0$, 
where $C>0$ is a constant independent of $\ep$. 
}

\vspace{5mm}
The paper is organized as follows. 
In Section 2 we prove the theorem for the case $\al=1$. 
The first step of our proof is to show the property of finite speed of propagation for classical solutions. We next show a generalized Kato's lemma to prove  Theorem 1.1. 
In Section 3 we consider the case $\al>1$. The theorem is proved in the same way as in Section 2. 
In Section 4 we apply the theorem to the original equation \eqref{ore}, and discuss the effect of the scale factor to the solutions.

\section{Case $\al=1$. }

\setcounter{equation}{0}

\quad 
We first consider the case $\al=1$. 
In order to prove Theorem 1.1 for $\al=1$,  we first show the following proposition: 

\noindent
{\bf Proposition 2.1. 
} 
{\it Let $F(u,u',u'')$ be a function of class $C^1$ in $u,u'$ and $u''$ satisfying 
\begin{equation}
F(0,0,u'')=0 \quad \mbox{for all }u''
\label{asm1}
\end{equation}
and let $u(t,x)$ be a $C^2$-solution of the equation
\begin{equation}
u_{tt}-t^{-2}\Delta u+\mu t^{-1} u_t=F(u,u',u''), 
\label{eqF1}
\end{equation}
in the region 
\[
\Lambda_{T,x_0}=\{(t,x)\in [1,T)\times \R^n: \; |x-x_0|<\ln T-\ln t\}
\]
for some $T>1$ and $x_0\in \R^n$. 
Assume that 
\begin{equation}
u(1,x)=u_t(1,x)=0 \qquad \mbox{for }|x-x_0|<\ln T. 
\label{asm2}
\end{equation}
Then $u$ vanishes in $\Lambda_{T,x_0}$. 
}
  
\pf  
We prove the proposition following \cite{Jo1,So}. 
Let 
\[
\psi(\lm,x)=\exp\lbt{\ln T-\lb{(\ln T-\lm)^2+(\ln T)^{-2}(2\lm\ln T-\lm^2)|x-x_0|^2}^{1/2}}. 
\]
We then have 
\[
\psi(0,x)=1, \quad \lim_{\lm\to \ln T}\psi(\lm,x)=\exp\{\ln T-|x-x_0|\}. 
\] 
Define the region $R_\lm$ by 
\[
R_\lm=\{(t,x): \; 1\le t\le \psi(\lm,x), \; |x-x_0|<\ln T-\ln t\}  
\]
and let $\lm_0$ satisfy $0<\lm_0<\ln T$. 
We have 
\[
\Lambda_{T,x_0}=\bigcup_{0\le \lm<\ln T}R_\lm
\]
and also 
\begin{equation}
|\nb_x\psi(\lm,x)|=\fa{(2\lm\ln T-\lm^2)\psi(\lm,x)|x-x_0|}{(\ln T)^2\lb{(\ln T-\lm)^2+(\ln T)^{-2}(2\lm\ln T-\lm^2)|x-x_0|^2}^{1/2}}
\label{psi}
\end{equation}
for $0\le \lm\le \lm_0<\ln T$.
Define the surface $S_\lm$ by  
\[
S_\lm=\{(t,x): t=\psi(\lm,x), \; |x-x_0|<\ln T\}. 
\]
We note that the outward unit normal at $(\psi(\lm,x),x)\in S_\lm$ is 
\[
\fa 1{\sr{1+|\nb_x\psi(\lm,x)|^2}}(1,-\nb_x\psi(\lm,x)). 
\]
Then 
\begin{align*}
\int_{R_\lm}2u_tFdtdx
&=\int_{R_\lm}2u_t(u_{tt}-t^{-2}\Delta u+\mu t^{-1}u_t)dtdx \\ 
&=\int_{R_\lm}\lb{\p_t(u_t^2+t^{-2}|\nb_xu|^2)-2\nb_x\cdot(t^{-2} u_t\nb_x u)
+2\al t^{-3}|\nb_x u|^2+2\mu t^{-1}u_t^2}dtdx \\
&\ge  \int_{R_\lm}\lb{\p_t(u_t^2+t^{-2}|\nb_xu|^2)-2\nb_x\cdot(t^{-2} u_t\nb_x u)}dtdx.  
\end{align*}
Note by \eqref{psi} that 
\[
|t^{-1}\nb_x\psi|<(\ln T)^{-1}\sr{2\lm_0\ln T-\lm_0^2}\equiv \theta(\lm_0)<1
\] 
on $S_\lm$ for $0\le \lm\le \lm_0<\ln T$. 
Using the divergence theorem, we obtain 
\begin{align}
\int_{R_\lm}2u_tFdtdx
&\ge  \int_{S_\lm}\lb{u_t^2+t^{-2}|\nb_xu|^2+2\nb_x\psi\cdot(t^{-2} u_t\nb_x u)}
\fa 1{\sr{1+|\nb_x\psi|^2}}d\sigma \nonumber \\
&\ge \int_{S_\lm}\lb{u_t^2+t^{-2}|\nb_xu|^2-t^{-1}|\nb_x\psi|(u_t^2+t^{-2}|\nb_x u|^2)}
\fa 1{\sr{1+|\nb_x\psi|^2}}d\sigma \nonumber \\
&\ge (1-\theta(\lm_0))\int_{S_\lm} \fa{u_t^2+t^{-2}|\nb_xu|^2}{\sr{1+|\nb_x\psi|^2}}d\sigma
\label{Ap1}
\end{align}
for $\lm\le \lm_0<\ln T$.


On the other hand, by assumption \eqref{asm1}, we have 
\[
|u_t F(u,u',u'')|\le C(u^2+|u'|^2)
\]
in $\Lambda_{T,x_0}$. 
We note that 
\[
\int_1^{\psi(\lm,x)}u(t,x)^2dt\le \fa 12\psi(\lm,x)^2\int_1^{\psi(\lm,x)}u_t(t,x)^2dt
\lsm T^2\int_1^{\psi(\lm,x)}u_t(t,x)^2dt
\]
and also that 
\begin{align*}
\psi_\lm(\lm,x)
&=\fa{\psi(\lm,x)(\ln T-\lm)(1-(\ln T)^{-2}|x-x_0|^2)}
{\lb{(\ln T-\lm)^2+(\ln T)^{-2}(2\lm\ln T-\lm^2)|x-x_0|^2}^{1/2}}, \\
|\psi_\lm(\lm,x)|
&\le C T. 
\end{align*}
We then have 
\begin{align}
\int_{R_\lm}2u_tFdtdx
&\le C(1+T^2)\int_{R_\lm} |u'|^2dxdt \nonumber\\
&\le C_T \int_{R_\lm} (u_t^2+t^{-2}|\nb_xu|^2)dxdt \nonumber\\
&=C_T\int_0^\lm \int_{S_\mu}(u_t^2+t^{-2}|\nb_xu|^2)\fa{\psi_\mu}{\sr{1+|\nb_x\psi|^2}}d\sigma d\mu \nonumber\\
&\le C_T\int_0^\lm \int_{S_\mu}\fa{u_t^2+t^{-2}|\nb_xu|^2}{\sr{1+|\nb_x\psi|^2}}d\sigma d\mu
\label{Ap2}
\end{align}
for $\lm \le \lm_0$. Set 
\[
I(\lm)=\int_{S_\lm} \fa{u_t^2+t^{-2}|\nb_xu|^2}{\sr{1+|\nb_x\psi|^2}}d\sigma. 
\]
Using Gronwall's inequality for \eqref{Ap1} and \eqref{Ap2}, we see that 
$I(\lm)=0$ for $0\le \lm\le \lm_0<\ln T$.  Therefore, since $\lm$ and $\lm_0$ are arbitrary, we see that $u'=0$  in $\Lambda_{T,x_0}$, and thus by \eqref{asm2} also that $u= 0$ in $\Lambda_{T,x_0}$. 
This completes the proof of the proposition. 
\hfill\qed 

\vspace{5mm}
From the proposition above, we easily see that the following corollary holds: 

\noindent
{\bf Corollary 2.2.} 
{\it 
Let $F$ in Proposition 2.1 satisfy \eqref{asm1}. 
If $u$ is a $C^2$-solution of \eqref{eqF1} and if 
$u(1,x)=u_t(1,x)=0$ for $|x|>R$ with $R>0$, then 
$u(t,x)=0$ for $|x|>R+\ln t$. 
}

\vspace{5mm}
We next show the following generalized Kato's lemma.

\noindent
{\bf Lemma 2.3.}
{\it Let $p>1, \;a\ge 0, \;b>0,\; c>0, \;q>0, \; \mu\ge 0$ and 
\[ 
M\equiv (p-1)(c-a)+2>0. 
\]
Let $T\ge T_1>T_0>1$. 
Assume that $F\in C^2([T_0,T))$ satisfies the following three conditions:  
\begin{align*}
(i) \quad & F(t) \ge A_0t^{-a}(\ln t)^{-b}(t-T_1)^c \qquad \mbox{for }t\ge T_1 \\
(ii) \quad &F''(t) +\fa{\mu F'(t)}{t}\ge A_1(\ln t)^{-q}|F(t)|^p \quad \mbox{for }t\ge T_0	 \\
(iii) \quad &F(T_0)\ge 0, \quad F'(T_0)>0,  
\end{align*}
where $A_0$ and $A_1$ are positive constants. Then 
$T$ has to satisfy
\[
T^{M/(p-1)}(\ln T)^{-b-q/(p-1)}<CA_0^{-1}, 
\]
where $C_0$ is a constant depending on $A_1,\mu,p,q,a,b$ and $c$. 
}
 
\pf We prove the lemma by referring to \cite{Jo2,TuLin1}.

Mutiplying assumption (ii) by $t^\mu$, we have  
\[
t^\mu F''+\mu t^{\mu-1}F'\ge A_1 t^\mu(\ln t)^{-q}|F|^p. 
\]
Integrating the above inequality over $[T_0,t]$ yields 
\begin{align}
t^\mu F'(t)-F'(T_0)\ge A_1\int_{T_0}^t s^\mu(\ln s)^{-q}|F(s)|^pds\ge 0. 
\label{e0}
\end{align}
By assumption (iii),  
\[
 F'(t)\ge t^{-\mu}F'(T_0)>0,
\]
thus $F'(t)>0$ for $t\ge T_0$. 
Moreover, 
\[
F'(t)\ge A_1 t^{-\mu}\int_{T_0}^t s^\mu(\ln s)^{-q}|F(s)|^pds, 
\]
hence, integrating again yields
\[
F(t)-F(T_0)  \ge A_1\int_{T_0}^t \tau^{-\mu}\int_{T_0}^\tau s^\mu(\ln s)^{-q}|F(s)|^pdsd\tau.  
\]
By assumptions (i) and (iii), we have
\begin{align*}
F(t) & \ge A_0^pA_1 \int_{T_1}^t \tau^{-\mu}\int_{T_1}^\tau s^{\mu-ap}(\ln s)^{-q-bp}(s-T_1)^{cp}dsd\tau \nonumber\\
& \ge A_0^pA_1 \int_{T_1}^t \tau^{-\mu-ap}(\ln\tau)^{-q-bp}\int_{T_1}^\tau (s-T_1)^{\mu+cp}dsd\tau \nonumber\\
&\ge \fa{A_0^pA_1}{\mu+cp+1}\int_{T_1}^t \tau^{-\mu-ap}(\ln\tau)^{-q-bp}(\tau-T_1)^{\mu+cp+1} ds 
\nonumber\\
 &\ge \fa{A_0^pA_1}{\mu+cp+1}t^{-\mu-ap}(\ln t)^{-q-bp}\int_{T_1}^t (\tau-T_1)^{\mu+cp+1} ds  \nonumber\\
 &= \fa{A_0^pA_1}{(\mu+cp+1)(\mu+cp+2)}t^{-\mu-ap}(\ln t)^{-q-bp}(t-T_1)^{\mu+cp+2} \nonumber\\
 &\ge \fa{A_0^pA_1}{(\mu+cp+2)^2}t^{-\mu-ap}(\ln t)^{-q-bp}(t-T_1)^{\mu+cp+2} \qquad \mbox{for }t\ge T_1. 
\end{align*} 
Based on the fact above, we define the sequences $a_j, \; b_j, \; c_j, \; D_j$ for $j = 0, 1, 2, \cdots$ by
\begin{align}
a_{j+1} = pa_j + \mu, & & b_{j+1} = pb_j+q,  
& & c_{j+1}=pc_j +\mu+2, & & D_{j+1} = \frac{A_1D_j^p}{(pc_j+\mu+2)^2}  \label{e2}\\
a_0 = a, & & b_0 = b, & & c_0 = c, & & D_0= A_0.  \label{e3}
\end{align}
Solving \eqref{e2} and \eqref{e3}, we obtain  
\begin{align*}
&a_j = p^j\lp{a+\fa{\mu}{p-1}}-\fa{\mu}{p-1}, \qquad 
b_j = p^j\lp{b+\fa{q}{p-1}}-\fa{q}{p-1},  \\
&c_j = p^j\lp{c+\fa{\mu+2}{p-1}}-\fa{\mu+2}{p-1}, 
\end{align*}
and thus
\[
D_{j+1} = \frac{A_1D_j^p}{c_{j+1}^2} \ge \lp{c+\fa{\mu+2}{p-1}}^{-2}\frac{A_1D_j^p}{p^{2j+2}}.
\]
Then, 
\begin{align*}
D_j &\ge \frac{BD_{j-1}^p}{p^{2j}} \\
&\ge \fa B{p^{2j}}\lp{\fa{BD_{j-2}^p}{p^{2(j-1)}}}^p =\fa{B^{1+p}}{p^{2j+2p(j-1)}}D_{j-2}^{p^2} \\
&\ge  \fa{B^{1+p}}{p^{2j+2p(j-1)}}\lp{\fa{BD_{j-3}^p}{p^{2(j-2)}}}^{p^2}
=\fa{B^{1+p+p^2}}{p^{2j+2p(j-1)+2p^2(j-2)}}D_{j-3}^{p^3} \\
&\ge \cdots\cdots \ge \fa{B^{1+p+p^2+\cdots+p^{j-1}}}{p^{2(j+p(j-1)+p^2(j-2)+\cdots+p^{j-1})}}D_0^{p^j}, \\
\intertext{and }
\ln D_j &\ge \ln B \sum_{k=0}^{j-1}p^k-2\ln p \sum_{k=0}^j kp^{j-k}+p^j\ln D_0 \\
 &= \frac{\ln B}{p-1}(p^j-1)-2p^j \sum_{k=0}^j 
\frac{k}{p^k}\ln p+p^j\ln D_0, 
\end{align*}
where $B=\{c+(\mu+2)/(p-1)\}^{-2}A_1$.  
For sufficiently large $j$, we have
\[
D_j \ge \exp(Ep^j),
\]
where
\begin{equation}
E = \frac{1}{p-1}\min\left(0, \ln B \right) - 2\sum_{k=0}^{\infty}\frac{k}{p^k}\ln p+\ln A_0. 
\label{eqofE}
\end{equation}
Thus, since $F(t)\ge D_jt^{-a_j}(\ln t)^{-b_j}(t-T_1)^{c_j}$ holds for $t\ge T_1$, we obtain
\begin{align}
F(t)  &\ge t^{\mu/(p-1)}(\ln t)^{q/(p-1)}(t-T_1)^{-(\mu+2)/(p-1)} \nonumber\\
  & \quad \cdot\exp\left[\left\{E+\lp{c+\fa{\mu+2}{p-1}}\ln(t-T_1)-\lp{a+\fa{\mu}{p-1}}\ln t-\lp{b+\fa{q}{p-1}}\ln(\ln t)\right\}p^j\right]
 \label{5}
\end{align}
for $t\ge T_1$. 
Since     
\[
\lp{c+\fa{\mu+2}{p-1}}-\lp{a+\fa{\mu}{p-1}}=c-a+\fa{2}{p-1}>0
\]
by assumption, choosing $t$ large enough, we can find a positive $\delta$ such that
\[
E+\lp{c+\fa{\mu+2}{p-1}}\ln(t-T_1)-\lp{a+\fa{\mu}{p-1}}\ln t-\lp{b+\fa{q}{p-1}}\ln(\ln t) \ge \delta > 0. 
\]
It then follows from \eqref{5}  that $F(t) \lr \infty$ as $j \to \infty$ for 
sufficiently large $t$. We therefore see that the life-span $T$ of $F(t)$ has to satisfy
\[
T^{M/(p-1)}(\ln T)^{-b-q/(p-1)}<CA_0^{-1},  
\]
where $M=(p-1)(c-a)+2>0$, and $C$ is a constant depending on $A_1,\mu,p,q,a,b$ and $c$. 
This completes the proof of the proposition. 
\hfill\qed

\vspace{2cm}
We are now in a position to prove Theorem 1.1 for $\al=1$ by applying Lemma 2.3. 
Set 
\[
F(t)=\int u(t,x)dx.
\] 
By Corollary 2.3, we have 
\[
\mbox{supp }u(t,\cdot)\subset \{|x|\le A(t)+R\}, 
\]
with
\[
A(t)=\int_1^t s^{-1}ds =\ln t. 
\]
By \eqref{eq}, 
\begin{align}
F''(t)+\fa\mu t F'(t) \ge \fa 1{(A(t)+R)^{n(p-1)}}|F(t)|^p 
\ge \fa C{(\ln t)^{n(p-1)}}|F(t)|^p. 
\label{Fp}
\end{align}
Mutiplying the inequality above by $t^\mu$ 
and integrating imply
\begin{equation}
t^\mu F'(t)-F'(1)\gtrsim \int_1^t \fa{s^{\mu}}{(\ln s)^{n(p-1)}}|F(s)|^pds\ge 0. 
\label{F01} 
\end{equation}
Note that 
\[
F'(t)>0 \quad \mbox{for }t\ge 1
\]
since $F'(1)>0$ by assumption. 
Moreover, from \eqref{F01}, 
\begin{align*}
F'(t)&\gtrsim t^{-\mu}F'(1). \\
\intertext{Integrating over $[1,t]$ implies}
F(t)-F(1)&\gtrsim F'(1)\int_1^t s^{-\mu}ds >0.
\end{align*}
We hence see that 
\[
F(t)\ge F(1)=C\ep>0 \qquad \mbox{for }t\ge 1
\]
by assumption. 
Going back to \eqref{F01}, we have  
\begin{align*}
t^\mu F'(t)-F'(1)&\ge C\ep^p (\ln t)^{-n(p-1)}\int_1^t (s-1)^\mu ds \\
 &\ge C\ep^p(\ln t)^{-n(p-1)}(t-1)^{\mu+1} \quad \mbox{for }t\ge 1. 
\end{align*}  
Since $F'(1)>0$, we have 
\[
F'(t)\ge C\ep^p (\ln t)^{-n(p-1)}t^{-\mu}(t-1)^{\mu+1} \quad \mbox{for }t\ge 1. 
\]
Integrating again, we have from $F(1)>0$ by assumption, 
\begin{align}
F(t)&\ge C\ep^p \int_1^t  (\ln s)^{-n(p-1)}s^{-\mu}(s-1)^{\mu+1}ds \nonumber \\
 &\ge C\ep^p (\ln t)^{-n(p-1)}t^{-\mu}\int_1^t (s-1)^{\mu+1}ds \nonumber \\
 &\ge C\ep^p (\ln t)^{-n(p-1)}t^{-\mu}(t-1)^{\mu+2}
 \qquad 
\mbox{for }t\ge 1. 
\label{Fup1} 
\end{align} 
Finally, by \eqref{Fp} and \eqref{Fup1}, applying Lemma 2.3 with $a=\mu$, $b=q=n(p-1)$, $c=\mu+2$ and $A_0=C\ep^p$, we obtain the desired results since
\[
M=2(p-1)+2=2p>0. 
\]
This completes the proof of Theorem 1.1 for $\al=1$. 
\hfill\qed

\section{
Case $\al>1$.}

\setcounter{equation}{0}
 
\quad 
We prove the theorem for the case $\al> 1$ in this section. 
We first show the following proposition: 

\noindent
{\bf Proposition 3.1. 
} 
{\it 
Let $F$ in Proposition 2.1 satisfy \eqref{asm1} 
and let $u(t,x)$ be a $C^2$-solution of the equation 
\begin{equation}
u_{tt}-t^{-2\al}\Delta u+\mu t^{-1} u_t=F(u,u',u''), 
\label{eqFa3l}
\end{equation}
with $\al>1$ 
in the region 
\[
\Lambda_{T,x_0}=\{(t,x)\in [1,T)\times \R^n: \; |x-x_0|<(t^{1-\al}-T^{1-\al})/(\al-1)\}
\]
for some $T>1$ and $x_0\in \R^n$. 
Assume that 
\begin{equation}
u(1,x)=u_t(1,x)=0 \qquad \mbox{for }|x-x_0|<\fa{1-T^{1-\al}}{\al-1}. 
\label{asm32}
\end{equation}
Then $u$ vanishes in $\Lambda_{T,x_0}$. 
}
  
\pf  
We prove the proposition above in a similar way to the proof of Proposition 2.1. 
Let 
\begin{align*}
&\psi(\lm,x)\\
=&\lbt{T^{1-\al}+\lb{(1-T^{1-\al}-\lm)^2+(1-T^{1-\al})^{-2}(2(1-T^{1-\al})\lm-\lm^2)(\al-1)^2|x-x_0|^2}^{1/2}}^{1/(1-\al)}. 
\end{align*}
We then have 
\[
\psi(0,x)=1, \quad \lim_{\lm\to 1-T^{1-\al}}\psi(\lm,x)=\{T^{1-\al}+(\al-1)|x-x_0|\}^{1/(1-\al)}. 
\]
Define the region $R_\lm$ by 
\[
R_\lm=\lb{(t,x): \; 1\le t\le \psi(\lm,x), \; |x-x_0|<(t^{1-\al}-T^{1-\al})/(\al-1)} 
\]
and let $0<\lm_0<1-T^{1-\al}$. 
We have 
\[
\Lambda_{T,x_0}=\bigcup_{0\le \lm<1-T^{1-\al}}R_\lm
\] 
and also 
\begin{align}
&|\nb_x\psi(\lm,x)| \nonumber\\
&=\fa{(1-T^{1-\al})^{-2}(2(1-T^{1-\al})\lm-\lm^2)\psi(\lm,x)^\al(\al-1)^2|x-x_0|}{(1-\al)\lb{(1-T^{1-\al}-\lm)^2+(1-T^{1-\al})^{-2}(2(1-T^{1-\al})\lm-\lm^2)(\al-1)^2|x-x_0|^2}^{1/2}}
\label{3psi}
\end{align}
for $0\le \lm\le \lm_0<1-T^{1-\al}$.
Define the surface $S_\lm$ by  
\[
S_\lm=\lb{(t,x): t=\psi(\lm,x), \; |x-x_0|<(1-T^{1-\al})/(\al-1)}. 
\]
Note by \eqref{3psi} that 
\[
|t^{-\al}\nb_x\psi|<(1-T^{1-\al})^{-1}\sr{2(1-T^{1-\al})\lm_0-\lm_0^2}\equiv \theta(\lm_0)<1
\] 
on $S_\lm$ for $0\le \lm\le \lm_0<1-T^{1-\al}$. 
Using the divergence theorem, we obtain 
\[
\int_{R_\lm}2u_tFdtdx
\ge (1-\theta(\lm_0))\int_{S_\lm} \fa{u_t^2+t^{-2\al}|\nb_xu|^2}{\sr{1+|\nb_x\psi|^2}}d\sigma
\]
for $\lm\le \lm_0<1-T^{1-\al}$.


On the other hand, since
\begin{align*}
\psi_\lm(\lm,x)
&=\fa{\psi(\lm,x)^\al(1-T^{1-\al}-\lm)(1-(1-T^{1-\al})^{-2}(\al-1)^2|x-x_0|^2)}
{(\al-1)\lb{(1-T^{1-\al}-\lm)^2+(1-T^{1-\al})^{-2}(2\lm(1-T^{1-\al})-\lm^2)(\al-1)^2|x-x_0|^2}^{1/2}}, \\
|\psi_\lm(\lm,x)|
&\le C T^\al, 
\end{align*}
we have 
\[
\int_{R_\lm}2u_tFdtdx
\le C_T\int_0^\lm \int_{S_\mu}\fa{u_t^2+t^{-2\al}|\nb_xu|^2}{\sr{1+|\nb_x\psi|^2}}d\sigma d\mu
\]
for $\lm \le \lm_0$. Therefore, we see that $u'=0$  in $\Lambda_{T,x_0}$, and thus by \eqref{asm2} also that $u= 0$ in $\Lambda_{T,x_0}$. 
This completes the proof of the proposition. 
\hfill\qed

\vspace{5mm}
From this proposition, we easily see that the following corollary holds: 

\noindent
{\bf Corollary 3.2.} 
{\it  
Let $F$ in Proposition 2.1 satisfy \eqref{asm1}. 
If $u$ is a $C^2$-solution of \eqref{eqFa3l} and if 
$u(1,x)=u_t(1,x)=0$ for $|x|>R$ with $R>0$, then 
$u(t,x)=0$ for $|x|>R+(1-t^{1-\al})/(\al-1)$. 
}

\vspace{5mm}
We next show another generalized Kato's lemma. 

\noindent
{\bf Lemma 3.3.}
{\it Let $p>1, \;a\ge 0, \;b>0,\; \mu\ge 0$ and 
\[ 
M\equiv (p-1)(b-a)+2>0. 
\]
Let $T\ge T_1>T_0>1$. 
Assume that $F\in C^2([T_0,T))$ satisfies the following three conditions:  
\begin{align*}
(i) \quad & F(t) \ge A_0t^{-a}(t-T_1)^b \qquad \mbox{for }t\ge T_1 \\
(ii) \quad &F''(t) +\fa{\mu F'(t)}{t}\ge A_1|F(t)|^p \quad \mbox{for }t\ge T_0	 \\
(iii) \quad &F(T_0)\ge 0, \quad F'(T_0)>0,  
\end{align*}
where $A_0$ and $A_1$ are positive constants. Then 
$T$ has to satisfy
\[
T<C_0A_0^{-(p-1)/M},  
\]
where $C_0$ is a constant depending on $A_1,\mu,p,a$ and $b$. 
}
 
\pf We prove the lemma above in a similar way to the proof of Lemma 2.3. 

Using assumption (ii), we have  
\[
F(t)-F(T_0)  \ge A_1\int_{T_0}^t \tau^{-\mu}\int_{T_0}^\tau s^\mu|F(s)|^pdsd\tau.  
\]
By assumptions (i) and (iii), we have
\begin{align*}
F(t) & \ge A_0^pA_1 \int_{T_1}^t \tau^{-\mu}\int_{T_1}^\tau s^{\mu-ap}(s-T_1)^{bp}dsd\tau \\
 &\ge \fa{A_0^pA_1}{(\mu+bp+2)^2}t^{-\mu-ap}(t-T_1)^{\mu+bp+2} \qquad \mbox{for }t\ge T_1. 
\end{align*} 
Proceeding as in the proof of Lemma 2.3, we obtain
\begin{align}
F(t)  &\ge t^{\mu/(p-1)}(t-T_1)^{-(\mu+2)/(p-1)} \nonumber\\
  & \quad \cdot\exp\left[\left\{E+\lp{b+\fa{\mu+2}{p-1}}\ln(t-T_1)-\lp{a+\fa{\mu}{p-1}}\ln t\right\}p^j\right]
 \label{5al}
\end{align}
for $t\ge T_1$, where $E$ is given in \eqref{eqofE} and $B=\{b+(\mu+2)/(p-1)\}^{-2}A_1$.    
Since     
\[
\lp{b+\fa{\mu+2}{p-1}}-\lp{a+\fa{\mu}{p-1}}=b-a+\fa{2}{p-1}>0
\]
by assumption, choosing $t$ large enough, we can find a positive $\delta$ such that
\[
E+\lp{b+\fa{\mu+2}{p-1}}\ln(t-T_1)-\lp{a+\fa{\mu}{p-1}}\ln t\ge \delta > 0. 
\]
It then follows from \eqref{5al}  that $F(t) \lr \infty$ as $j \to \infty$ for 
sufficiently large $t$. We therefore see that the life-span $T$ of $F(t)$ has to satisfy
\[
T<CA_0^{-(p-1)/M},  
\]
where $M=(p-1)(b-a)+2>0$, and $C$ is a constant depending on $A_1,\mu,p,a$ and $b$. 
This completes the proof of the lemma. 
\hfill\qed

\vspace{2cm}
We now prove Theorem 1.1 for $\al>1$ by applying Lemma 3.3. 
Set 
\[
F(t)=\int u(t,x)dx.
\] 
By Corollary 3.2, we have 
\[
\mbox{supp }u(t,\cdot)\subset \{|x|\le A(t)+R\}, 
\]
with
\[
A(t)=\int_1^t s^{-\al}ds =\fa 1{\al-1}(1-t^{1-\al}). 
\]
By \eqref{eq}, 
\begin{align}
F''(t)+\fa\mu t F'(t) \ge \fa 1{(A(t)+R)^{n(p-1)}}|F(t)|^p 
\ge C|F(t)|^p. 
\label{3Fp}
\end{align}
Proceeding as in the proof for the case $\al =1$, we obtain
\begin{align*}
t^\mu F'(t)-F'(1)
&\ge C \int_1^t s^{\mu}|F(s)|^pds \\
&\ge C\ep^p \int_1^t (s-1)^\mu ds \\
 &\ge C\ep^p(t-1)^{\mu+1} \quad \mbox{for }t\ge 1. 
\end{align*}  
Since $F'(1)>0$, we have 
\[
F'(t)\ge C\ep^p t^{-\mu}(t-1)^{\mu+1} \quad \mbox{for }t\ge 1. 
\]
Integrating again, we have from $F(1)>0$ by assumption, 
\begin{align}
F(t)&\ge C\ep^p \int_1^t  s^{-\mu}(s-1)^{\mu+1}ds \nonumber \\
 &\ge C\ep^p t^{-\mu}\int_1^t (s-1)^{\mu+1}ds \nonumber \\
 &\ge C\ep^p t^{-\mu}(t-1)^{\mu+2}
 \qquad 
\mbox{for }t\ge 1. 
\label{3Fup1} 
\end{align} 
Finally, by \eqref{3Fp} and \eqref{3Fup1}, applying Lemma 3.3 with $a=\mu$, $b=\mu+2$ and $A_0=C\ep^p$, we obtain the desired results since
\[
M=(p-1)(b-a)+2=2p>0. 
\]
This completes the proof of Theorem 1.1 for $\al>1$. 
\hfill\qed

\newpage

\section{Wave Equation in FLRW}
\setcounter{equation}{0}

\quad 
In this section we summarize the results obtained so far for the original equation \eqref{ore}, which corresponds to \eqref{eq} with $\al=2/(n(1+w))$ and $\mu=2/(1+w)$. 

For the case of decelerated expansion $2/n-1<w\le 1, \; n\ge 2$, we have shown in \cite{TW1,TW2} that blow-up in finite time occurs if $1<p\le p_F(n-2/(1+w))=1+2/\{n-2/(1+w)\}$ or if $1<p\le p_c(n,w)$. Moreover, the following upper bounds of the lifespan are proved: 
\begin{align}
&T_\ep\le C\ep^{\fa{-(p-1)}{2-(n-2/(1+w))(p-1)}}  & &  \mbox{if } 1<p<p_F(n-2/(1+w)), 
\label{31} \\
&T_\ep\le \exp(C\ep^{-(p-1)}) & & \mbox{if } p=p_F(n-2/(1+w)), 
\nonumber \\
&T_\ep\le C\ep^{\fa{-2p(p-1)}{\gamma_0(n,p,w)}}  &&  \mbox{if } 1<p<p_c(n,w),  
\label{32}\\
&T_\ep\le \exp(C\ep^{-p(p-1)}) &&  \mbox{if } p=p_c(n,w)>p_F(n-2/(1+w)), \nonumber
\end{align}
where $p_c(n,w)$ is the positive root of the equation $\gamma_0(n,p,w)=0$, and 
\begin{align*}
\gamma_0(n,p,w)&=\lp{1-\fa 2{n(1+w)}}\gamma\lp{n,p,\fa 2{n(1+w)},\fa 2{1+w}}\\
&=-(n-1)p^2+\lp{n+1+\fa 4{n(1+w)}}p+2-\fa 4{n(1+w)}. 
\end{align*}
Recall that $\gamma(n,p,\al,\mu)$ is given in \eqref{gm}.  
Note that $p_c(n,w)>p_S(n)$ and $\gamma_0(n,p,w)>\gamma_S(n,p)
$, where $p_S(n)$ is the Strauss exponent and $\gamma_S(n,p)$ is given in \eqref{snp}.

On the other hand, for  $-1<w \le 2/n-1$ and $n\ge 2$, which represents an accelerated expanding universe, we can apply Theorem 1.1 to the original equation \eqref{ore} since $\al\ge  1$ and $\mu\ge n$ in \eqref{eq}.  
By those theorems, we see that blow-up in finite time can happen to occur for all $p>1$. This is in contrast to the case of decelerated expansion above. We also recall that the lifespan satisfies
\begin{equation}
\begin{cases}
T_\ep^2(\ln T_\ep)^{-n(p-1)}\le C\ep^{-(p-1)} &\mbox{if }w = 2/n-1, \\
T_\ep\le C\ep^{-(p-1)/2} & \mbox{if }-1<w < 2/n-1, 
\end{cases}
\label{33}
\end{equation}
where $C>0$ is a constant independent of $\ep$.  
  
Figures 1, 2 and 3 show the range of blow-up conditions in terms of $w$ and $p$ in the case $n=2,3,5$, respectively. 
Note that if $p_F(n-2/(1+w))=p_c(n,w)$, then $w$ is the larger root $w^\ast$ of the equation
\[
n(n^2+n+2)w^2+2n(n-1)^2w+n^3-5n^2+8n-8=0. 
\]
Each Region (A) is for the case of the accelerated expanding universe, where the life span of blow-up solutions is dominated by \eqref{33}. 
In contrast, each Region (B) and (C) represents the decelerated expanding one. 
In Region (B) \eqref{32} is better than \eqref{31}, on the other hand, this relation becomes reverse in Region (C).

From these results, we see that the blow-up range of $p$ in the flat FLRW spacetime is larger than that in the Minkowski spacetime. 
Moreover, as is the case in the decelerated expanding, the lifespan of the blow-up solutions in the the accelerated expanding universe is shorter than that in the Minkowski spacetime since
$\ep^{-(p-1)/2}<\ep^{-2p(p-1)/\gamma(n,p,0,0)}$ for sufficiently small $\ep$. Therefore, we can say that finite time blow-up  can occur more easily in the FLRW spacetime. 

We will treat the remaining case $w=-1$ in future papers.  
\begin{figure}[H]
\includegraphics[width=13cm, bb=00 400 500 760, clip]{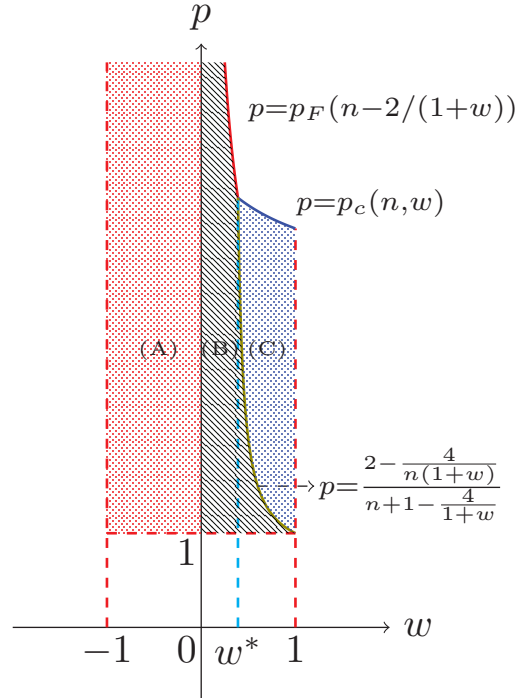}
\caption{Range of blow-up conditions in case $n=2$}

\end{figure}

\begin{figure}[p]
\includegraphics[width=13cm, bb=00 500 500 800, clip]{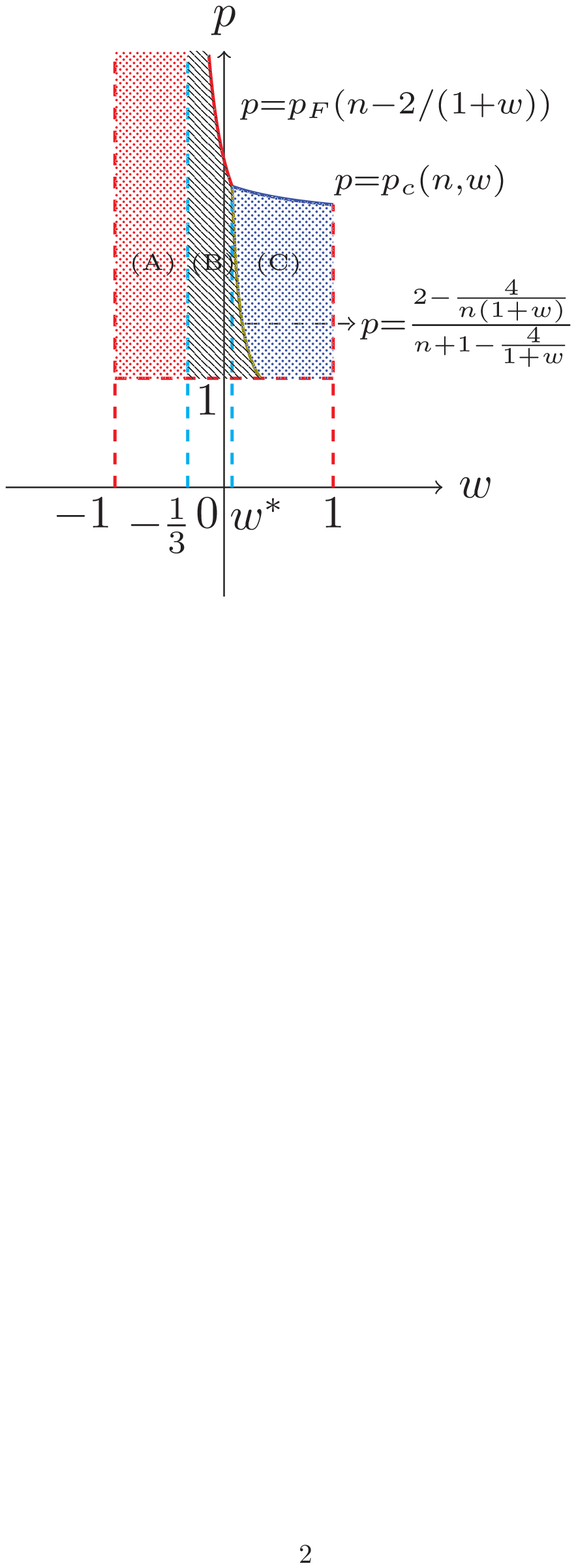}
\caption{Range of blow-up conditions in case $n=3$}

\includegraphics[width=15cm, bb=00 500 500 800, clip]{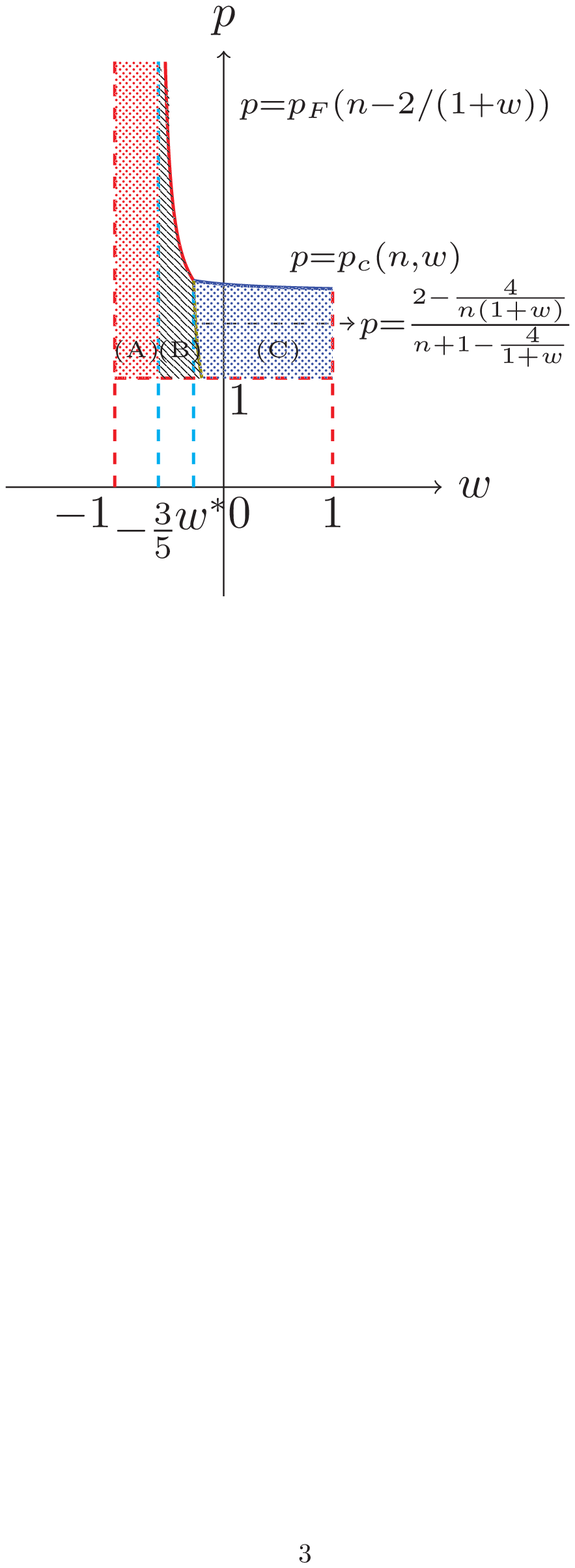}
\caption{Range of blow-up conditions in case $n=5$}

\end{figure}

\newpage


\end{document}